\documentclass[12pt]{article}

\usepackage[latin1]{inputenc}
\usepackage[english]{babel}
\usepackage{amssymb,amsmath,amscd,amstext,amsthm}
\usepackage{amsfonts}
\usepackage{amsthm}
\usepackage[dvips]{graphicx}
\usepackage{psfrag}
\usepackage{epsfig}
\usepackage{palatino}
\usepackage{xypic}
\input{xy}
\xyoption{all}

\newcounter{eno}
\newcounter{ch}
\newcounter{sec}
\newcounter{no}
\def \sec#1 {\stepcounter{sec} \bigskip\par{\underline{\bf \S \arabic{sec} {\bf #1}  \setcounter{eno}{1} \setcounter{no}{1}}}}
\def \prop#1 {\noindent{\bf Proposition \arabic{sec}.\arabic{no}} \stepcounter{no} {\it #1}}
\def \cor#1 {\noindent{\bf Corollary \arabic{sec}.\arabic{no}} \stepcounter{no} {\it #1}}
\def \lem#1 {\noindent{\bf Lemma \arabic{sec}.\arabic{no}} \stepcounter{no} {\it #1}}
\def \thm#1 {\noindent{\bf Theorem \arabic{sec}.\arabic{no}} \stepcounter{no} {\it #1}}
\def \conj#1 {\noindent{\bf Conjecture \arabic{sec}.\arabic{no}} \stepcounter{no} {\it #1}}
\def \define#1 {\noindent{\bf Definition \arabic{sec}.\arabic{no} } \stepcounter{no} {\it #1} \hspace{1mm}}

\def\co{\colon\thinspace}

\newcommand {\Z} {\mathbb Z}

\newcommand {\ZG} {\mathbb {Z}[G]}

\newcommand {\ZD} {\mathbb {Z}[D_{4n}]}

\newcommand {\Zd} {\mathbb {Z}[D_{8}]}

\def \proof#1 { \hspace{2mm}Proof:$\,\,\,$  #1  ${}$\hfill $\Box$ \,\, \newline }

\begin{document}

\begin{center}
{\huge \bf The D(2) property for $D_8$}

{W.H.Mannan}

  \begin{quote}Wall's D(2) problem asks if a cohomologically 2- dimensional geometric 3- complex is necessarily homotopy equivalent to a geometric 2- complex.  We solve 
  part of the problem when the fundamental group is dihedral of order $2^n$, and offer a complete solution for the case where it is $D_8$ the dihedral group of order 8.

\end{quote}
  
  \bigskip
  {\bf MSC} {57M20, 57M05} \hfill  {\bf Keywords} {Wall's D(2) problem, algebraic }
  
  \hfill{complexes, k-- invariants}

\end{center}

%Section 1
\sec{Introduction}

\hspace{2mm}Wall introduces the D(2) problem in \cite{Wall}.  The D(2) property has been verified for dihedral groups of order $4n+2$ in \cite{John}.  Hence we consider dihedral 
groups of 
order $4n$.  As these do not have periodic resolutions, not all the methods of \cite{John} can be applied to them.  Our main result is orthogonal to the result in \cite{John1}, in
the sense that
it concerns dihedral groups whose order is a power of $2$, rather than twice an odd number. 

\hspace{2mm}  We begin by recalling some of the theory of $k$-- invariants.  We work over a finite group $G$, of order $n$.

\define{Algebraic complex} We define an algebraic $n$-- complex, to be a sequence of maps and modules:

$$
F_n \stackrel{d_n}{\to} F_{n-1} \stackrel{d_{n-1}}{\to} \cdots \to F_1 \stackrel{d_1}{\to} F_0
$$

where the $F_i$ are free finitely generated modules over $\ZG$, the cokernel of $d_1$ is $\Z$ (with trivial $G$ action), and the sequence is exact at $F_1$.

\hspace{2mm}Let $(F_i, d_i)$  and $(F_i',d_i')$, $i=0,1,2$, denote algebraic 2-- complexes.  Let the $f_i \co F_i\to F_i'$, $i=0,1,2$ constitute a chain map, $f$, between them.

\prop{{\rm (See Prop. 47.1 of \cite{John1})} $f$ is a homotopy equivalence if and only if it induces isomorphisms ker$(d_2) \to $ ker$(d_2')$ and coker$(d_1)\to$ coker$(d_1')$.}

\define{Algebraic $\pi_2$} We define $\pi_2(F_i, d_i)$ to be ker$(d_2)$.

\hspace{2mm}Let $J$ denote the kernel of $d_2$ and let $J'$ denote the kernel of $d_2'$.

\prop{(i) Given $\alpha \co J \to J'$, we may choose a chain map \newline $f_\alpha \co (F_i, d_i) \to (F_i',d_i')$ which induces $\alpha \co J \to J'$. A map $\Z\to\Z$, is 
induced on the cokernels.  We
denote this map by multiplication by $k$.
\newline
(ii)   The congruence of $k$ modulo $n$ is independent of the choice of $f_\alpha$.  
\newline
(iii) Given $k'$ congruent to $k$ modulo $n$, we may choose a chain map $f'_\alpha$, which also induces $\alpha \co J \to J'$, and which induces multiplication by $k'$ on $\Z$.}

\hspace{2mm}Proof:  i) See Prop 25.3 of \cite{John1}.

\hspace{12mm}  ii) {\rm See Prop 25.3 of \cite{John1} and  see Prop 33.3 of \cite{John1}}.

\hspace{12mm}  iii) Let  $\epsilon \co
F_0 \to F_0 /{\rm Im}(d_i) \cong \Z$, $\epsilon'\co  F_0' \to F_0' /{\rm Im}(d_i') \cong \Z$ denote the natural quotient maps.  Pick $x \in F_0'$ such that $\epsilon'x=1$.  Let
$h \co \Z \to F_0'$ be the map sending $1 \in \Z$ to $\sum_{g \in G}xg$.  Then 

$$\epsilon'h(1)=\epsilon'(\sum_{g \in G}xg)= \sum_{g\in G} (\epsilon'x) g = \sum_{g\in G}1=n$$

\hspace{2mm}Let $(f_\alpha')_2 =(f_\alpha)_2$

\hspace{7mm} $(f_\alpha')_1 =(f_\alpha)_1$
 
\hspace{7mm}  $(f_\alpha')_0 =(f_\alpha)_0 + (\frac{k'-k}{n}) h\epsilon$.  
  
\hspace{2mm}  Then $f_\alpha'$ is a chain map as 

$$(f_\alpha')_0 d_1 = (f_\alpha)_0 d_1 +  (\frac{k'-k}{n}) h\epsilon d_1 = (f_\alpha)_0 d_1 +0 = d_1' (f_\alpha)_1 =d_1' (f_\alpha')_1$$

\hspace{2mm}Finally note 

$$\epsilon'(f_\alpha')_0 = \epsilon'(f_\alpha)_0  +  \epsilon'(\frac{k'-k}{n}) h\epsilon=k\epsilon + (k'-k)\epsilon = k' \epsilon$$

\hfill $\Box$

\define{k-- invariant} Given $\alpha$ as in the proposition, we define the congruence of $k$ modulo $n$ to be $k_\alpha$.  

\hspace{2mm}We have a ring homomorphism $\kappa \co $End$(J)\to \Z_n$, sending $\alpha \mapsto k_\alpha$.

\lem{{\rm (See Prop. 26.6 of \cite{John1})} The kernel of $\kappa$ is all maps which factor through a projective module.}

\lem{{\rm (See Prop. 33.7 of \cite{John1})}$\kappa$ is independent of the choice of algebraic complex $(F_i,d_i)$.}

\proof{$\kappa_1=1$, so $\kappa$ is surjective.  Hence $\kappa$ is equal to the quotient map End$(J) \to$ End$(J)/$Ker$(\kappa)$ composed 
with a ring isomorphism $\Z_n \to \Z_n$.  However, any ring isomorphism $\Z_n \to \Z_n$
must map $1 \mapsto 1$.  Hence it must be the identity.}

\define{Swan map} The Swan map is the homomorphism $Aut(J) \to \Z_n^*$ which sends an automorphism to its k-- invariant.

\prop{If the Swan map Aut$(J) \to \Z_n^*$ is surjective, and we have an isomorphism $\alpha \co J \to J'$, then $(F_i, d_i)$ ard $(F_i', d_i')$ are chain homotopy equivalent.}

\proof{By surjectivity we may choose $\beta \co J \to J$, such that $k_\beta=k_\alpha^{-1}$.  Then by proposition 1.4(iii), we may pick $f_{\alpha\beta}$ which induces isomorphisms
$J \to J'$ and the identity $\Z \to \Z$.  Hence by proposition 1.2, $f_{\alpha\beta}$ is a homotopy equivalence.}

\lem{Given a map $\alpha \co J \to J$, let $\alpha' \co J \oplus\ZG \to J \oplus\ZG$ denote the map $\alpha \oplus 1$.  Then $k_\alpha =k_{\alpha'}$.}

\hspace{2mm}  Hence it is sufficient to show that the Swan map is surjective for $J$, in order to deduce that it is surjective for $J \oplus \ZG^r$, for all natural numbers $r$.
 Consequently we have: 
 
\prop{ If the Swan map is surjective for $J$, then there is a unique algebraic $2$-- complex, up to chain homotopy 
equivalence, with algebraic $\pi_2$ equal to $J\oplus \ZG^r$, for each $r$.}

 \hspace{2mm}For $n$ coprime to $3$, we will show that in the case $G= \ZD$, the unit 
 \newline
 $3\in \Z_{4n}$ is in the image of the Swan map for $J$, where $J$ is the algebraic $\pi_2$
 of a particular algebraic $2$-- complex.  We will then show that $-1$ and $3$ generate the units of $\Z_{2^n}$, so the Swan map is surjective for $J$, for dihedral groups of order
 $2^n$. So there is a unique algebraic 2-- complex, up to chain homotopy equivalence, with algebraic $\pi_2$ equal to $J\oplus \ZD^r$, for each $r$. 
 
 \hspace{2mm}  We then show that $J$ has minimal $\Z$-- rank, for a module which occurs as an algebraic $\pi_2$.  We use a cancellation result due to Swan \cite{Swan2} to show
 that for the group $D_8$, the only modules which arise as an algebraic $\pi_2$ of an algebraic $2$-- complex, are of the form $J \oplus \ZD^r$. We have shown by 
 this point, that there is only
 one algebraic $2$-- complex which has each of these algebraic $\pi_2$'s, up to chain homotopy equivalence.  We show each of these are geometrically realized.  
 
 \hspace{2mm}  Finally, we quote Theorem I of \cite{John1} which states that if every algebraic $2$-- complex over a finite group $G$ is geometrically realized, then $G$
 satisfies the D(2) property.

%Section 2
\sec{Surjectivity of the Swan map.}

\hspace{2mm} Let $D_{4n}$ be the group given by the presentation, $\langle a,b \,\, |a^{2n}=b^2=e,\,\,aba =b \rangle$.  $\Sigma$ will denote 
$\sum_{i=0}^{2n-1} a^i$.  This presentation has a Cayley complex, which in turn has an associated algebraic complex.  This is an exact
sequence over $\ZD$:

$$
J\hookrightarrow \ZD^3 \stackrel{\partial_2}{\longrightarrow} \ZD^2 \stackrel{\partial_1}{\longrightarrow} \ZD \stackrel{\epsilon}{\twoheadrightarrow} \Z
\eqno (1)
$$

$\epsilon$ is determined by mapping $1 \in \ZD$ to $1 \in \Z$.  $J$ is the kernel of $\partial_2$. Let $e_1$, $e_2$ denote 
basis elements of $\ZD^2$.  Then $\partial_1 e_1 = a-1, \, \, \partial_1 e_2 = b-1$.  

\hspace{2mm} Let $E_1, \, E_2, \, E_3$ be basis elements of $\ZD^3$, which correspond to the relations in the presentation so that:

$
\partial_2 E_1 = e_1 \Sigma \newline
\partial_2 E_2 = e_2 (1+b) \newline
\partial_2 E_3 = e_1 + e_2a + e_1ba - e_2= e_1(1+ba)+e_2(a-1)
$

\hspace{2mm}  With respect to the basis $\{E_1,\, E_2,\,E_3\}$ and the basis $\{e_1, \, e_2\}$, $\partial_2$ is given by ;

$$
\left[ \begin{array}{ccc} \Sigma & 0 & 1+ba \\ 0 & 1+b & a-1 \end{array} \right]
$$

Let

$\alpha_0 = 1+a+b$

$\alpha_1 = \left[ \begin{array}{cc} 1+a-ba & b-1 \\ 0 & 1 \end{array} \right]$

$\alpha_2 = \left[ \begin{array}{ccc}  1+a-ba & 0&0 \\ 0 & 1&0 \\ 0&0&1\end{array} \right]$

\hspace{2mm}It is easily verified that:

\prop{The following diagram commutes: 
\newline
\begin {eqnarray*}
J\hookrightarrow \ZD^3 \stackrel{\partial_2}{\longrightarrow} \ZD^2 \stackrel{\partial_1}{\longrightarrow}
\ZD \stackrel{\epsilon}{\twoheadrightarrow} \Z \,\,\,\\
\downarrow \theta  \,\,\,\,\, \quad \downarrow \alpha_2 \qquad \qquad \downarrow \alpha_1 \qquad \downarrow \alpha_0 \qquad\downarrow 3
\\
J\hookrightarrow \ZD^3 \stackrel{\partial_2}{\longrightarrow} \ZD^2 \stackrel{\partial_1}{\longrightarrow}
\ZD \stackrel{\epsilon}{\twoheadrightarrow} \Z \,\,\,
\end {eqnarray*}
\newline
where $\theta$ is the restriction of $\alpha_2$.}

\hspace{2mm} For the remainder we will assume $3$ coprime to $n$.  Our goal is to show that $\theta$ is an isomorphism. As we know that $\kappa_\theta=3$, this will suffice to
show that $3$ is in the image of the Swan map.  

\hspace{2mm}Note that if we regard the above diagram as a diagram of commutative $\Z$--
modules and $\Z$- linear maps,  there are well defined integer determinants for all the maps in the chain map.  A map is an 
isomorphism if and only if it has
determinant $\pm1$.  (As the property of being an isomorphism is dependent only on surjectivity and injectivity, it does not depend on
whether we are regarding modules as being over $\ZD$, or $\Z$). 

\hspace{2mm} Note also that over $\Z$, all the maps in the exact sequences above are quotienting of a summand, followed by inclusion
of a summand.  Consequently, the following proposition holds;

\prop{$3{\rm Det}(\theta) {\rm Det}(\alpha_1)= {\rm Det}(\alpha_2){\rm Det}(\alpha_0)$}

\hspace{2mm} Proof:  Let $u$ be the restriction of $\alpha_1$ to the kernel of $\partial_1$ and let $v$ be the restriction of
$\alpha_0$ to the kernel of $\epsilon$.    Then by the previous discussion, we have 

$$
3{\rm Det}(\theta){\rm Det}(\alpha_1)= {\rm Det}(\theta){\rm Det}(u) {\rm Det}(v)3 = {\rm Det}(\alpha_2){\rm Det}(\alpha_0)
$$
\hfill $\Box$

\hspace{2mm} We will use this to show that $Det(\theta)=1$.

\prop{${\rm Det}(1+a+b)=-3$}

\hspace{2mm} Proof:  Let $A$ be the matrix for left multiplication by $1+a+b$ in the regular representation, with basis
$\{a^{2n-1},a^{2n-2},\cdots, a,1,ba^{2n-1},ba^{2n-2},\cdots, ba,b\}$.  Then the upper right quadrant of $A$ and the lower left quadrant
of $A$ are copies of the identity matrix.  The upper left quadrant has $1$'s along the diagonal and immediately
above as well as a $1$ in the bottom left corner.  The lower right quadrant has $1$'s along the diagonal and
immediately below, as well as a $1$ in the top right corner.  All the other entries in $A$ are $0$.

\hspace{2mm} For example, if $n$ were equal to $4$, the matrix $A$ would be

$$ \tiny
\left[\begin{array}{ccccccccccccccccc}  
1&1&0&0&0&0&0&0&&1&0&0&0&0&0&0&0\\
0&1&1&0&0&0&0&0&&0&1&0&0&0&0&0&0 \\
0&0&1&1&0&0&0&0&&0&0&1&0&0&0&0&0\\
0&0&0&1&1&0&0&0&&0&0&0&1&0&0&0&0 \\
0&0&0&0&1&1&0&0&&0&0&0&0&1&0&0&0\\
0&0&0&0&0&1&1&0&&0&0&0&0&0&1&0&0 \\
0&0&0&0&0&0&1&1&&0&0&0&0&0&0&1&0\\
1&0&0&0&0&0&0&1&&0&0&0&0&0&0&0&1 \\
 & & & & & & & && & & & & & & & \\
1&0&0&0&0&0&0&0&&1&0&0&0&0&0&0&1\\
0&1&0&0&0&0&0&0&&1&1&0&0&0&0&0&0 \\
0&0&1&0&0&0&0&0&&0&1&1&0&0&0&0&0\\
0&0&0&1&0&0&0&0&&0&0&1&1&0&0&0&0 \\
0&0&0&0&1&0&0&0&&0&0&0&1&1&0&0&0\\
0&0&0&0&0&1&0&0&&0&0&0&0&1&1&0&0 \\
0&0&0&0&0&0&1&0&&0&0&0&0&0&1&1&0\\
0&0&0&0&0&0&0&1&&0&0&0&0&0&0&1&1 \\
\end{array} \right]
$$

\hspace{2mm} Label the rows of $A$, $v_1, v_2.....v_{4n}$. We will perform row
operations.

\hspace{2mm}  First let $v_{2n}'= v_{2n}-v_1+v_2-v_3....-v_{2n-1}$.  Now let $v_{2n}''=v_{4n}$ and $v_{4n}''=v_{2n}'$.  Let the
remaining $v_{i}''=v_{i}$.  This swap causes a change of sign in the determinant, so the matrix with rows $v_i''$ has determinant -Det$A$. 
In the case $n=4$, this matrix is

$$ \tiny
\left[\begin{array}{ccccccccccccccccc}  
1&1&0&0&0&0&0&0&&1&0&0&0&0&0&0&0\\
0&1&1&0&0&0&0&0&&0&1&0&0&0&0&0&0 \\
0&0&1&1&0&0&0&0&&0&0&1&0&0&0&0&0\\
0&0&0&1&1&0&0&0&&0&0&0&1&0&0&0&0 \\
0&0&0&0&1&1&0&0&&0&0&0&0&1&0&0&0\\
0&0&0&0&0&1&1&0&&0&0&0&0&0&1&0&0 \\
0&0&0&0&0&0&1&1&&0&0&0&0&0&0&1&0\\
0&0&0&0&0&0&0&1&&0&0&0&0&0&0&1&1 \\
 & & & & & & & && & & & & & & & \\
1&0&0&0&0&0&0&0&&1&0&0&0&0&0&0&1\\
0&1&0&0&0&0&0&0&&1&1&0&0&0&0&0&0 \\
0&0&1&0&0&0&0&0&&0&1&1&0&0&0&0&0\\
0&0&0&1&0&0&0&0&&0&0&1&1&0&0&0&0 \\
0&0&0&0&1&0&0&0&&0&0&0&1&1&0&0&0\\
0&0&0&0&0&1&0&0&&0&0&0&0&1&1&0&0 \\
0&0&0&0&0&0&1&0&&0&0&0&0&0&1&1&0\\
0&0&0&0&0&0&0&0&&-1&1&-1&1&-1&1&-1&1 \\
\end{array} \right]
$$

\hspace{2mm}  For each $2n+1 \leq i \leq 4n-2$, let $v_i'''=v_i''+v_{i+1}''-v_{i-2n}''$.  

\hspace{2mm} Let $v_{4n-1}'''=
v_{4n-1}''+v_{2n}''-v_{2n-1}$''.  

\hspace{2mm}For $i \leq 2n$ let $v_i'''=v_i''$.

\hspace{2mm}When $n=4$, the matrix with rows $v_i'''$ is

$$ \tiny
\left[\begin{array}{ccccccccccccccccc}  
1&1&0&0&0&0&0&0&&1&0&0&0&0&0&0&0\\
0&1&1&0&0&0&0&0&&0&1&0&0&0&0&0&0 \\
0&0&1&1&0&0&0&0&&0&0&1&0&0&0&0&0\\
0&0&0&1&1&0&0&0&&0&0&0&1&0&0&0&0 \\
0&0&0&0&1&1&0&0&&0&0&0&0&1&0&0&0\\
0&0&0&0&0&1&1&0&&0&0&0&0&0&1&0&0 \\
0&0&0&0&0&0&1&1&&0&0&0&0&0&0&1&0\\
0&0&0&0&0&0&0&1&&0&0&0&0&0&0&1&1 \\
 & & & & & & & && & & & & & & & \\
0&0&0&0&0&0&0&0&&1&1&0&0&0&0&0&1\\
0&0&0&0&0&0&0&0&&1&1&1&0&0&0&0&0 \\
0&0&0&0&0&0&0&0&&0&1&1&1&0&0&0&0\\
0&0&0&0&0&0&0&0&&0&0&1&1&1&0&0&0 \\
0&0&0&0&0&0&0&0&&0&0&0&1&1&1&0&0\\
0&0&0&0&0&0&0&0&&0&0&0&0&1&1&1&0 \\
0&0&0&0&0&0&0&0&&0&0&0&0&0&1&1&1\\
0&0&0&0&0&0&0&0&&-1&1&-1&1&-1&1&-1&1 \\
\end{array} \right]
$$

\hspace{2mm}
In general, the matrix with rows $v_i'''$ has an upper
triangular top left quadrant, with $1$'s along the diagonal and a lower left quadrant with no non-zero entries.  Let $B$ denote the
lower right quadrant.  Then ${\rm Det}(1+a+b)=-{\rm Det}(B)$.     

\hspace{2mm} Cycle the top $2n-1$ rows of $B$ upwards to get the matrix $B'$.  As this is a cycle of odd length, ${\rm Det}(B')={\rm
Det}(B)$.  When $n=4$, $B'$ is

$$ \tiny
\left[\begin{array}{ccccccccc}  
1&1&1&0&0&0&0&0 \\
0&1&1&1&0&0&0&0\\
0&0&1&1&1&0&0&0 \\
0&0&0&1&1&1&0&0\\
0&0&0&0&1&1&1&0 \\
0&0&0&0&0&1&1&1\\
1&1&0&0&0&0&0&1\\
-1&1&-1&1&-1&1&-1&1 \\
\end{array} \right]
$$

\hspace{2mm} Label the rows of $B'$ as $w_1,...,w_{2n}$.  Set $u_{i}=w_i-w_{i+1}$ for $i=1,2, \cdots,2n-3$.  Let $B''$ denote the matrix with
rows $u_i$.  After these row operations, we have ${\rm Det}(1+a+b)=-{\rm Det}(B'')$

\hspace{2mm} When $n=4$, $B''$ is

$$ \tiny
\left[\begin{array}{ccccccccc}  
1&0&0&-1&0&0&0&0 \\
0&1&0&0&-1&0&0&0\\
0&0&1&0&0&-1&0&0 \\
0&0&0&1&0&0&-1&0\\
0&0&0&0&1&0&0&-1 \\
0&0&0&0&0&1&1&1\\
1&1&0&0&0&0&0&1\\
-1&1&-1&1&-1&1&-1&1 \\
\end{array} \right]
$$

\hspace{2mm}  We must consider two cases: $n$ congruent to $1$ modulo $3$ and $n$ congruent to $2$ modulo $3$.

\hspace{2mm}  If $n=1$ modulo $3$ then replace $u_{2n-1}$ with 

$u_{2n-1}-u_1-u_2-u_4-u_5-u_7-u_8\cdots - u_{2n-3}$.

\hspace{2mm}Also, replace $u_{2n}$ with 

$u_{2n} + (u_1-u_2+u_3)+(u_7-u_8+u_9)+(u_{13}-u_{14}+u_{15})...+(u_{2n-7}-u_{2n-6}+u_{2n-5})$.  

We are left with a
matrix with $1$'s along the diagonal and $0$'s below, except in the last four columns.  The 4 by 4 matrix in the bottom right corner is 

$$\left[\begin{array}{cccc} 1&0&0&-1 \\ 0&1&1&1 \\ 0&0&1&2 \\ 0&0&-1&1 \end{array} \right] $$

$${\rm Det}(1+a+b)=-{\rm Det}\left[\begin{array}{cccc} 1&0&0&-1 \\ 0&1&1&1 \\ 0&0&1&2 \\ 0&0&-1&1 \end{array} \right] 
=-{{\rm Det}}\left[\begin{array}{cccc} 1&0&0&-1 \\ 0&1&1&1 \\ 0&0&1&2 \\ 0&0&0&3 \end{array} \right]=-3$$

\hspace{2mm}  If $n=2$ modulo $3$ then replace $u_{2n-1}$ with 

$u_{2n-1}-u_1-u_2-u_4-u_5-u_7-u_8....-u_{2n-5}$.  

\hspace{2mm}Also, replace $u_{2n}$ with 

$u_{2n} + (u_1-u_2+u_3)+(u_7-u_8+u_9)+(u_{13}-u_{14}+u_{15})...+(u_{2n-9}-u_{2n-8}+u_{2n-7})$.  

We are left with a
matrix with $1$'s along the diagonal and $0$'s below, except in the last four columns.  The 4 by 4 matrix in the bottom right corner is 

$$\left[\begin{array}{cccc} 1&0&0&-1 \\ 0&1&1&1 \\ 1&1&0&1 \\ -1&1&-1&1 \end{array} \right] $$

$${\rm Det}(1+a+b)=-{\rm Det}\left[\begin{array}{cccc} 1&0&0&-1 \\ 0&1&1&1 \\ 1&1&0&1 \\ -1&1&-1&1 \end{array} \right] 
=-{\rm Det}\left[\begin{array}{cccc} 1&0&0&-1 \\ 0&1&1&1 \\ 0&1&0&2 \\ 0&1&-1&0 \end{array} \right]$$

$$=-{\rm Det}\left[\begin{array}{cccc} 1&0&0&-1 \\ 0&1&1&1 \\ 0&0&-1&1 \\ 0&0&-2&-1 \end{array} \right] = 
-{\rm Det}\left[\begin{array}{cccc} 1&0&0&-1 \\ 0&1&1&1 \\ 0&0&-1&1 \\ 0&0&0&-3 \end{array} \right] = -3$$

\hfill $\Box$

\prop{${\rm Det}(2-b)=3^{2n}$}

\hspace{2mm} Proof:  Let $A$ be the matrix for $2-b$ in the regular representation, with basis
$\{1, b, a,ba, a^2, ba^2......a^{2n-1},ba^{2n-1} \}$ .  Then $A$ consists of $2n$ two by two blocks of the form

$$\left[\begin{array}{cc}  2 & -1 \\ -1 & 2 \end{array} \right]$$

along the diagonal.  Hence ${\rm Det} (A) = 3^{2n}$

\hfill $\Box$

\prop{${\rm Det}(1+a-ba) \neq 0$}

\hspace{2mm} Proof:  Let $\alpha_2'=\left[\begin{array}{ccc}  2-b & 0&0 \\ 0 & 1&0 \\ 0&0&1\end{array}
\right]$

The following diagram commutes:

\begin {eqnarray*}
J\hookrightarrow \ZD^3 \stackrel{\partial_2}{\longrightarrow} \ZD^2 \stackrel{\partial_1}{\longrightarrow}
\ZD \stackrel{\epsilon}{\twoheadrightarrow} \Z \,\,\,\\
\downarrow \eta  \,\,\,\, \quad \downarrow \alpha_2' \qquad \qquad \downarrow \alpha_1 \qquad \downarrow \alpha_0 \qquad\downarrow 3
\\
J\hookrightarrow \ZD^3 \stackrel{\partial_2}{\longrightarrow} \ZD^2 \stackrel{\partial_1}{\longrightarrow}
\ZD \stackrel{\epsilon}{\twoheadrightarrow} \Z \,\,\,
\end {eqnarray*}

where $\eta$ is the restriction of $\alpha_2'$

Therefore $3{\rm Det}(\eta) {\rm Det}(\alpha_1)= {\rm Det}(\alpha_2'){\rm Det}(\alpha_0)$

So $3*{\rm Det}(\eta){\rm Det}(1+a-ba) = -3*3^{2n}$.  Hence ${\rm Det}(1+a-ba)$ cannot be $0$.
\hfill $\Box$

\prop{$\theta$ is an isomorphism}

\hspace{2mm} Proof: $3{\rm Det}(\theta) {\rm Det}(\alpha_1)= {\rm Det}(\alpha_2){\rm Det}(\alpha_0)$

\hspace{2mm} Therefore $3{\rm Det}(\theta){\rm Det}(1+a-ba)= -3 {\rm Det}(1+a-ba)$.  As ${\rm Det}(1+a-ba)$ is non-zero, we can conclude that

\hspace{2mm} ${\rm Det}(\theta)=-1$.

\hspace{2mm}Hence $\theta$ is an isomorphism.

\hfill $\Box$ $\,\,\,$

\cor{If $3 \in (\Z_{4n})^*$ then $3$ is in the image of the Swan Map:\newline ${\rm Aut(J)} \to (\Z_{4n})^*$ }

\hspace{2mm}  Let us now consider dihedral groups of order $2^m$ for $m \geq 2$.  Clearly $2^m$ is divisible by $4$
and coprime to $3$.  Hence we know that $3$ is in the image of the Swan Map.

\lem{$2^m$ divides $3^{2^{m-3}}-1+2^{m-1}$ for $m \geq 4$.}

\hspace{2mm} Proof:  We proceed by induction.  $3^{2^{4-3}}-1+2^{4-1} = 16$.  So the proposition holds for $m=4$. 
Now suppose it holds for some $m$.  Then 

$2^mz =3^{2^{m-3}}-1+2^{m-1}$ for some $z$.  Rearranging gives $3^{2^{m-3}}=1-2^{m-1}+2^mz $.  Then squaring gives: 

$$
3^{2^{m+1-3}}= (3^{2^{m-3}})^2 = (2^mz+1-2^{m-1})^2
$$

So

$$
3^{2^{m+1-3}}-1+2^{m+1-1} = (2^mz+1-2^{m-1})^2 -1 + 2^m
$$

$$
= 2^{2m}z^2+2^{2m-2} +2^{m+1}z-2^{2m}z = 2^{m+1}(2^{m-1}(z^2-z) + 2^{m-3} + z)
$$

So the proposition holds for $m+1$.  Hence by induction it holds for all $m\geq 4$.
\hfill $\Box$

\prop{The elements $3$, $-1$ generate $(\Z/{2^m})^*$ for $m \geq 2$.}

\hspace{2mm} Proof:  The order of $(\Z/{2^m})^*$ is $2^{m-1}$.  $(\Z/{4})^* = \{1,3\}$ and $(\Z/{8})^* = \{1, -1,3,
-3\}$, so only the case $m \geq 4$ remains.  We know that the order of $3$ in $(\Z/{2^m})^*$ is a power of $2$.  The
previous lemma shows us that for $m \geq 4$ it is at least $2^{m-2}$, as

$3^{2^{m-3}} \equiv 1+2^{m-1}$ Mod $2^m$.

\hspace{2mm}It remains to show that $-1$ is not a power of $3$, as
then the $\pm 3^k$ give us all $2^{m-1}$ elements of $(\Z/{2^m})^*$.

\hspace{2mm} Suppose $3^k= -1 \quad {\rm Mod} \,\, 2^m$ for some $m \geq 4$.  Then $3^k = -1 \quad {\rm Mod} \,\, 8$ which is impossible as $3^k$ only
takes the values $1$ and $3$ modulo $8$.

\hfill $\Box$

\hspace{2mm}Combining this result with corollary 2.7 we obtain

\cor{The Swan Map $Aut(J) \to (\Z_{2^m})^*$ is surjective for all $m \geq 2$.}

\hspace{2mm}From proposition 1.11, we may conclude

\thm{Over $\Z[D_{2^m}]$ an algebraic $2$-- complex, $X$,with \newline $\pi_2(X)=J\oplus \Z[D_{2^m}]^r$  is unique up to chain homotopy equivalence.}

%Section 3
\sec{The D(2) property for $\Zd$} 

\hspace{2mm}  Let $\mathbb{F}_2$ denote the two element module over $\ZD$, on which the
action of $\ZD$ is trivial.

\prop{({\rm See \cite{Adem}, p127}) (i)$H^0(D_{4n},\mathbb{F}_2)=\mathbb{F}_2$
\newline ${}$ \hspace{51mm}(ii)$H^1(D_{4n},\mathbb{F}_2)=\mathbb{F}_2^2$
\newline ${}$ \hspace{51mm}(iii)$H^2(D_{4n},\mathbb{F}_2)=\mathbb{F}_2^3$}

\hspace{2mm}Recall the sequence(1), from \S2.  By Schanuel's lemma, any module occurring as the algebraic $\pi_2$ of an algebraic $2$-- complex, over $\ZD$, must be stably
equivalent to $J$.

\prop{$J$ has minimal $Z$-- rank in its stable class.}

\hspace{2mm} Proof:  Given any finite algebraic 2-- complex, consider the cochain obtained by applying ${\rm Hom}_{Z[D_{4n}]} (\bullet
,\mathbb{F}_2)$:

$$
\mathbb{F}_2^{d_2} \stackrel{v_2}{\leftarrow} \mathbb{F}_2^{d_1} \stackrel{v_1}{\leftarrow} \mathbb{F}_2^{d_0} 
$$

where ${d_0}$, ${d_1}$, ${d_2}$, are the $\ZD$ ranks of the modules
in the complex.  As $H^0(D_{4n},\mathbb{F}_2)=\mathbb{F}_2$, the kernel
of $v_1$ has $\mathbb{F}_2$-- rank $1$.  Consequently, the image of  $v_1$ has $\mathbb{F}_2$-- rank
$d_0-1$.  $H^1(D_{4n},\mathbb{F}_2)=\mathbb{F}_2^2$ so $v_2$ has kernel of $\mathbb{F}_2$-- 
rank $2+d_0-1=d_0+1$.  The image of $v_2$ is then seen to have rank $d_1-d_0-1$.  
$H^2(D_{4n},\mathbb{F}_2)=\mathbb{F}_2^3$ so we know that $d_2 \geq 3+ d_1- d_0 -1$. 
Rearranging gives $d_2-d_1+d_0 \geq 2$.  

\hspace{2mm} Exactness implies that the $\Z$-- rank of the algebraic $\pi_2$ of the algebraic complex must be $4n(d_2-d_1+d_0)-1$. Hence our inequality implies that this is at 
least $8n-1$, which is the $\Z$-- rank of $J$.  

\hfill $\Box$

\hspace{2mm}We now restrict to the case $n=2$.

\prop{The only elements in the stable class of $J$ are modules
of the form $J \oplus \Zd^k$.}

\hspace{2mm}  Proof:  We refer to [3], Theorem 6.1.  This states that over $\Zd$, $A \oplus C =B \oplus C$ implies $A=B$ for torsion
free, finitely generated modules $A$, $B$, $C$.  

\hspace{2mm} If a module $M$ is in the stable class of $J$ then $M \oplus \Zd^r=J \oplus \Zd^s$.  From 
proposition 3.2 we have $s \geq r$.  From the theorem, we deduce that $M= J \oplus \Zd^{s-r}$.
\hfill $\Box$
\newline

\thm{The group $D_8$ satisfies the D(2) property.}

\hspace{2mm} Proof:  The only modules that can turn up as the algebraic $\pi_2$ of an algebraic 2-- complex over $\Zd$ are ones of the form $J \oplus
\Zd^s$ for some $s \geq 0$.  Theorem 2.11 tells us that for 
each $s$, up to chain homotopy equivalence, there is a unique algebraic 2-- complex with algebraic $\pi_2$ equal to $J \oplus
\Zd^s$.  Given any $r$, the chain
homotopy equivalence class of this algebraic $2$-- complex is realized by the Cayley complex of the presentation:

$$
\langle a,b \,\, |a^{2n}=b^2=e,\,\,aba =b, \,\, r_1=e,\, r_2=e, \dots r_s=e \rangle
$$

where $r_i = e$ for $i=1, \dots ,s$.  

\hspace{2mm}Hence we know that every algebraic $2$-- complex over $D_8$ is geometrically realized.  By Theorem I of \cite{John1}, this is equivalent to $D_8$ satisfying the D(2)
property.

\hfill $\Box$

  \bibliographystyle{gtart}

 \noindent Address:University College London 
 \newline
 Gower Street 
 \newline 
  London WC1E 6BT
\newline  
email: \verb|wajid@mannan.info|
  \end{document}